\documentclass[12pt]{article}
\usepackage{amsmath}
\usepackage{amsfonts}
\usepackage{amssymb}
\usepackage[mathscr]{eucal}
\newcounter{theorem}
\newcounter{lemma}
\newcounter{definition}
\newcounter{remark}
\newcounter{corollary}
\newcounter{proposition}

\newenvironment{lemma}{%
      \refstepcounter{lemma}{\noindent\textsf{LEMMA\  \thelemma}.}%
      \bgroup\slshape}{\egroup}
\newenvironment{definition}{%
     \refstepcounter{definition}%
     {\noindent\textsf{DEFINITION\ \thedefinition.}}%
      \bgroup\slshape}{\egroup}

\newenvironment{corollary}{%
      \refstepcounter{corollary}{\noindent\textsf{COROLLARY\ \thecorollary}.}%
      \bgroup\slshape}{\egroup}

\def\thetheorem{\arabic{section}.\arabic{theorem}}
\def\thelemma{\arabic{section}.\arabic{lemma}}
\def\thedefinition{\arabic{section}.\arabic{definition}}
\def\theremark{\arabic{section}.\arabic{remark}}
\def\thecorollary{\arabic{section}.\arabic{corollary}}
\def\theproposition{\arabic{section}.\arabic{proposition}}
\makeatletter
\renewcommand{\section}{\@startsection{section}{1}{0.0mm}%
{-\baselineskip}{0.5\baselineskip}{\normalfont\Large\bfseries}}
\makeatother \DeclareMathOperator{\Vol}{\textup{Vol}}
\begin{document}
\title{\textbf{Steiner-Minkowski Polynomials of Convex Sets
in High Dimension,\\
and Limit Entire Functions.}}
\author{\textbf{Victor Katsnelson}}
\date{ }
\maketitle
\begin{abstract}
For a convex set \(K\) of the \(n\)-dimensional Euclidean space,
the Steiner-Minkowski polynomial \(M_K(t)\) is defined as the
\(n\)-dimensional Euclidean volume of the neighborhood of the
radius \(t\). Being defined for positive \(t\), the
Steiner-Minkowski polynomials are considered for all complex
\(t\). The renormalization procedure for Steiner polynomial is
proposed. The renormalized Steiner-Minkowski polynomials
corresponding to all possible solid convex sets in all dimensions
form a normal family in the whole complex plane. For each of the
four families of convex sets: the Euclidean balls, the cubes, the
regular cross-polytopes and the regular symplexes of dimensions
\(n\), the limiting entire functions  are calculated explicitly .
\end{abstract}

\section{Renormalized Steiner-Minkowski \\Polynomials.}
\setcounter{equation}{0}%
\noindent%
NOTATION. \(\kappa_l\) is the \(n\)-dimensional volume of the
unite ball in \(\mathbb{R}^n\):
\begin{equation}%
\kappa_l=\frac{\pi^{l/2}}{\Gamma(\frac{l}{2}+1)}
\end{equation}
Let \(K\), \(K\subset\mathbb{R}^n\), be a compact convex set. For
\(t>0\), there defined the function
\begin{equation}
\label{ODMP}
M_K(t)\stackrel{\textup{def}}{=}\Vol_n(K+tB^n),\quad t>0.
\end{equation}
It is known (H.Minkowski) that \(M_K(t)\), considered for \(t>0\),
is a polynomial in \(t\) of degree \(n\):
\begin{equation}
\label{CMP}%
 M_K(t)=\sum\limits_{0\leq{}l\leq{}n}m_j(K)t^l\,.
\end{equation}
 This polynomial, which
was defined originally by \eqref{ODMP} for \(t>0\) only, will be
considered for all \textit{complex} \(t\).\\

\begin{definition}
The polynomial \(M_K(t)\) is said to be \textsf{the Steiner-Minkowski polynomial}
of the set \(K\).
\end{definition}

The following normalizations for the coefficients of the
polynomial \(M_K(t)\) are common:
\begin{equation}
\label{CSM}%
M_K(t)=\sum\limits_{0\leq{}l\leq{}n}\binom{n}{l}W_l(K)t^l\,,
\end{equation}
and
\begin{equation}
\label{IV}%
M_K(t)=\sum\limits_{0\leq{}l\leq{}n}\kappa_lV_{n-l}(K)\,t^l\,.
\end{equation}
So,
\begin{equation}
\label{VN}%
m_l(K)=\binom{n}{l}W_l(K)=\kappa_lV_{n-l}(K),\quad
0\leq{}l\leq{}n\,.
\end{equation}

The value \(W_j(K)\) is said to be \textit{l-th
cross-sectional measure} (or \textit{l-th
quermassintegral} - in the German manner) of the set \(K\).

The value \(V_l(K)\) is said to be \textit{l-th intrinsic
volume} of the set \(K\).

The constant term \(m_0(K)\) of the polynomial \(M_K(t)\) is the
\(n\)-dimensional volume of \(K\), the coefficient \(m_1(K)\) of
its linear term is the \((n-1)\)-dimensional `area' of the
`boundary surface' \(\partial{}K\):
\begin{subequations}
\label{FC}
\begin{alignat}{6}
m_0(K)&=\Vol_n(K), & \quad&(=&&W_0(K)&=&&&V_n(K)&&),
\label{FC0}\\
 m_1(K)&=\Vol_{n-1}(\partial{}K)&\quad &(=n&&W_1(K)&=&\kappa_1&&V_{n-1}(K)&&).
 \label{FC1}
\end{alignat}
\end{subequations}

We introduce the following normalization of the Steiner-Minkowski
polynomial of a convex set \(K\), under the \textsf{extra assumption}:\\
The set \(K\) is \textsf{solid}, that is the interior \(K\) is non-empty.
This extraassumption can be reformulated as:
\begin{equation}
\label{EA}%
\Vol_n(K)>0\,.
\end{equation}
Under the extra assumption \eqref{EA}, the area of the surface is
automatically positive: \(\Vol_{n-1}(\partial{}K)>0\,.\)\\

 The ratio
\begin{equation}
\label{SF}%
 \sigma_K=\frac{\Vol_{n-1}(\partial{}K)}{\Vol_n(K)}
\end{equation}
is said to be the \textit{shape factor} of the set \(K\). In terms
of cross-sectional measures the shape-factor is expressed as
\begin{equation}
\label{SFa}%
\sigma_K=\frac{n\,W_1(K)}{W_0(K)}\,.
\end{equation}
The shape factor has dimension \((\textup{length})^{-1}\).

\begin{definition}
Given a solid compact convex set \(K\) , we
present the Steiner-Minkovski polynomial \(M_K(t)\) in the form:
\begin{equation}
\label{DNMP} M_K(t)=\Vol_n(K)\cdot\mathscr{M}_K(\tau),
\end{equation}
where
\begin{equation}
\label{Rnt}%
\tau=\sigma_K\,t
\end{equation}
is a dimensionless parameter, and \(\mathscr{M}_K(\tau)\) is a
polynomial in \(\tau\) of degree \(n\).

The polynomial \(\mathscr{M}_K(\tau)\) in the variable \(\tau\) is
said to be \textsf{the renormalized Steiner-Minkowski polynomial
for the set \(K\).}
\end{definition}

Let us present the sequence of the coefficients of the normalizes
Steiner-Minkowski polynomial in the form:
\begin{equation}
\label{NSMP}%
 \mathscr{M}_K(\tau)=
 \sum\limits_{0\leq{}l\leq{}n}j_{n,l}\,\frac{\mu_l(K)}{l!}\,\tau^l\,,
\end{equation}
where the factors \(j_{n,l}\), so called \textit{Jensen
multipliers}, are \footnote{In other words,
\(\displaystyle\frac{j_{n,\,l}\,n^l}{l!}=\binom{n}{l}\) for
\(0\leq{}l\leq{}n\).}
\begin{equation}%
\label{JF}%
j_{n,0}=1,\quad j_{n,l}=\prod\limits_{0\leq{r}\leq{l-1}}\left(1-\frac{r}{n}\right),\ \ \text{for}
\ \ l=1,2,\,\dots\,,\,n,\quad j_{n,l}=0 \ \ \text{for}
\ \ l>n\,.
\end{equation}%
The Jensen multipliers possesses properties
\begin{equation}
\label{PJM}%
 0\leq{}j_{n,\,l}\leq{}1, \ \ \text{for all}\
n,\,l,\quad{}j_{n,\,l}\to{}1,\  \text{as \(l\) is fixed,}\ \
n\to\infty\,.
\end{equation}
\begin{definition}
The coefficients \(\mu_l(K)\) which appear in \eqref{NSMP} are said to be the \textsf{renormalized
Steiner-Minkowski coefficients} for the convex set \(K\).\\
\end{definition}

In view of \eqref{CMP},\eqref{FC}, \eqref{SFa}, \eqref{DNMP},
\eqref{Rnt},
\begin{subequations}
\label{ESC}
\begin{equation}
\label{ESC.a}
\mu_l(K)=\frac{(W_0(K))^{l-1}W_l(K)}{(W_1(K))^l}\,,\quad 0\leq{}l\leq{}n\,.
\end{equation}
In particular,
\begin{equation}
\label{ESC.b}
\mu_0(K)=1,\quad \mu_1(K)=1\,.
\end{equation}
\end{subequations}

\begin{lemma} For any solid compact convex set \(K\), \(K\subset\mathbb{R}^n,\)
the sequence \(\mu_l(K),\ l=0,\,1,\,2,\,\dots,\,n\)
of its \textsf{renormalized
Steiner-Minkowski coefficients} possesses the properties:
\begin{enumerate} \item
\addtocounter{equation}{1}
\begin{math}%
\label{FB}
\hfill 0<\mu_l(K)\leq{}1\,,\quad 0\leq l\leq n\,.\hfill (\arabic{section}.\arabic{equation})
\end{math}
\item
It is \textup{logarithmically concave}, that is the inequalities
\begin{equation}
\label{RAFI}
(\mu_l(K))^2\geq{}\mu_{l-1}(K)\mu_{l+1}(K),\quad 1\leq{}l\leq{}n-1
\end{equation}
hold.
\end{enumerate}
\end{lemma}

\noindent%
\textsf{PROOF.}  For any convex set \(K,\,K\subset{}\mathbb{R}^n,\)
 its cross-sectional measures
\( W_l(K),\,0\leq{}l\leq{}n,\) are non-negative, and if the interior of \(K\)
is non-empty, they are strictly positive:
\begin{equation}
\label{Po}
     W_l(K)>0,\quad{}0\leq{}l\leq{}n.
\end{equation}
The positivity property \eqref{Po} is a special case of the positivity
property of mixed volumes. (See, for example, \cite{}.) The inequalities
\(0<\mu_l(K)\,,\quad 0\leq l\leq n\,,\)  are   consequence of \eqref{ESC.a}
and  \eqref{Po}.

The cross-sectional measures of any convex set satisfy the Alexandrov-Fenchel
inequalities
\begin{equation}
\label{AFIn}%
W_j^2(K)\geq{}W_{j-1}(K)\,W_{j+1}(K),\ \ 1\leq j\leq n-1\,.
\end{equation}
The inequalities  \eqref{RAFI} are the Alexandrov-Fenchel inequalities \eqref{AFIn},
recalculated, according to \eqref{ESC.a}.

The inequalities \(\mu_l(K)\leq{}1\,,\ 0\leq l\leq n\,,\)   are consequence
 of the equalities \eqref{ESC.b} and  the logarithmic concavity inequalities  \eqref{RAFI}.\\

The following statement is an immediate consequence of \eqref{NSMP}  and  \eqref{FB}:

\vspace*{1.0ex}
\begin{lemma}
The renormalized Stener-Minkowski polynomial \(\mathscr{M}_K(\tau)\) of any solid compact convex set \(K\)
 satisfies the inequality
\begin{equation}
    |\mathscr{M}_K(\tau)|\leq{}\exp{\{|\tau|\}}  \quad \textup{for every} \ \ \tau\in\mathbb{C}\,.
\end{equation}
\end{lemma}

\begin{corollary}
The family \(\{\mathscr{M}_K(\tau)\}\)  of all renormalized Steiner Minkowski
 polynomials corresponding to all  solid compact convex sets \(K\subset\mathbb{R}^n\)
 in any dimension \(n:1\leq{}n<\infty\), is a normal family of analytic functions in
 the complex plane \(\mathbb{C}\).
\end{corollary}

This property of normality  prompts us the following statement of the problem:

\textit{Given a sequence \(\{K_p\}_p\)  of solid compact sets,
\(K_p\subset\mathbb{R}^{n_p}\), of increasing dimension:
 \(n_1<n_2<n_3\,\dots\) \,.  It is required
to study the cluster set of functions for the family \(\{\mathscr{M}_{K_p}(\tau)\}\)
of the corresponding Steiner-Minkowski polynomials.
  In particular, one need to clarify
under which conditions this cluster set consists of one function, that is
there exists the limit}
\begin{equation}
      \mathscr{M}(\tau) =\lim_{p\to\infty}\mathscr{M}_{K_p}(\tau)\,.
\end{equation}

Every function \(\mathscr{M}(\tau)\) from this cluster set is an entire function
satisfying the condition
\begin{equation}
|\mathscr{M}(\tau)|\leq \exp\{|\tau|\} \quad \textup{for all}\ \ \tau\in\mathbb{C}\,.
\end{equation}

It is natural to restrict our consideration to the case
\begin{equation}
       K_1\subset{}K_2\subset\,\cdots \,\subset{}K_p\,\cdots \, ,
\end{equation}
where we may consider the ambient Euclidean spaces as naturally  embedded:
\begin{equation}
    \mathbb{R}^{n_1}\subset \mathbb{R}^{n_2}\subset\,\cdots \,\subset\mathbb{R}^{n_p}\subset\,\cdots \,.
\end{equation}

 \section{Some Examples.}
 \setcounter{equation}{0}
    To calculate explicitly the Steiner-Minkowski polynomials for
 concrete families of convex sets is difficult.  Here we take several examples
 where such calculations can be done.\\

\noindent%
\textbf{1. The family of Euclidean balls \(B^n\).}
\begin{equation}
\label{KisB}
B^n=\{x=(\xi_1,\,\xi_2,\,\dots\,,\xi_n)\in\mathbb{R}^{n}:
\sum\limits_{1\leq{l}\leq{}n}|\xi_l|^2\leq{}\rho^2\}\,,\ \ \rho>0.
\end{equation}
Since \(B^n+tB^n=(\rho+t)B_n\) for \(t>0\),
\[\Vol_n(B^n+tB^n)=\Vol_n(B^n)(\rho+t)^n\quad \textup{for}\ \ t>0\,.\]
Thus the Steiner-Minkowski polynomial \(M_{B^n}(t)\) is
\[M_{B^n}(t)=\kappa_n(\rho+t)^n\,.\]
The shape constant \(\sigma_{B_n}\) for the \(n\)-dimensional ball
of radius \(\rho\) is:
\begin{equation}
\label{SCB}%
 \sigma_{B^n}=n/\rho\,.
\end{equation}
Thus, the renormalized Steiner-Minkowski polynomial
\(\mathscr{M}_{B^n}(\tau)\) is equal to:
\begin{equation}
\label{RMPB}%
\mathscr{M}_{B^n}(\tau)=(1+\tau/n)^n\,.
\end{equation}
The limit %
 \(\displaystyle\mathscr{M}_{B^\infty}(\tau) \stackrel{\text{\tiny def}}{=}
 \lim_{n\to\infty}\mathscr{M}_{B^n}(\tau)\) is:
\begin{equation}
\label{RMPBI}%
\mathscr{M}_{B^\infty}(\tau)=\exp\{\tau\}\,,
\end{equation}
or, in the term of Taylor series,
\begin{equation}
\label{RMPBIo}%
\mathscr{M}_{B^\infty}(\tau)=\sum\limits_{0\leq{}l<\infty}\frac{1}{\Gamma(l+1)}\,\tau^l\,,
\end{equation}

\centerline{\(\diamond\)\quad\(\diamond\)\quad\(\diamond\)}

The next three families of convex sets which we consider are the
families of cubes, regular crosspolytopes and regular symplexes. All these sets
are regular polytopes. The Minkowski polynomial for the polytope \(K,\,K\subset{}\mathbb{R}^n,\)
can be expessed in terms of the intrinsic volumes
\(V_r(K)\), \(r=0,\,1,\,\dots\,,n\), by \eqref{IV}.
 The \(r\)-th intrinsic volume \(V_r(K)\) of the
polytope \(K\) can be calculated by the formula
\begin{equation}
\label{FIV} V_r(K)=\sum\limits_{F_r}\gamma(F_r)\Vol_r(F_r),
\end{equation}
where the sum is taken over all \(r\)-faces \(F_r\) of
the polytope \(K\), and \(\gamma(F_r)\) is the \textit{external angle} at
the face \(F_r\), normalized so that the total angle is \(1\).
(See \cite{Gr}, Chapter 14, or \cite{Schn}). For regular polytope \(K\),
the \(r\)-volumes \(\Vol_r(F_r)\) of all its \(r\)-faces \(F_r\) are equal,
and all external angles \(\gamma_{F_r}\) are equal.
Their common value is denoted by \(v_r\) and \(\gamma_r\) respectively:
\begin{gather}
\label{vr}
\Vol_r(F_r)=v_r\quad \ \ \textup{for any}\ \ \ r-\text{face} \ F_r\,,\\
\gamma_(F_r)=\gamma_r\quad  \ \ \textup{for any}\ \ \ r-\text{face} \ F_r\,.
\label{gr}
\end{gather}
The cardinality of the set \(\mathscr{F}_r\) of all \(r\)-faces is denoted by \(\nu_r\):
\begin{equation}
\label{Car}
\nu_r=\#(\mathscr{F}_r).
\end{equation}
For the \textit{regular} polytope \(K\), the formula \eqref{FIV} takes the form
\begin{equation}
\label{IVRP}
V_r(K)=\nu_r\gamma_r v_r\,.
\end{equation}
So, the formula \eqref{IV} takes the form
\begin{equation}
\label{MPRP}
M_K(t)=\sum\limits_{0\leq{}l\leq{}n}\kappa_l\,\nu_{n-l}\,\gamma_{n-l}\,v_{n-l}\,t^l\,.
\end{equation}
In the dimension \(n\), \(n\geq{}5\), there are only three regular polytopes: cube,
crosspolytpe and symplex. We calculate the Minkowski polynomials for every of these
three families. It is easy to calculate the volumes and the total numbers of \(r\)-faces,
however to calculate the external angles for regular crosspolytopes and symplexes
is more difficult. These values can not be expressed `in elementary functions'.
In theyr expressions the Gauss' integral oshibok appears.

\centerline{\(\diamond\)\quad\(\diamond\)\quad\(\diamond\)}

\noindent%
 \textbf{2. The family of cubes \(Q^n\).}
\begin{equation}
\label{KisQ}
Q^n=\{x=(\xi_1,\,\xi_2,\,\dots\,,\xi_n)\in\mathbb{R}^{n}:
|\xi_l|\leq{}\rho,\ \ 1\leq{}l\leq{}n\}\,,\ \ \rho>0\,.
\end{equation}
The total number of \(l\)-faces \(\nu_l(Q^n)\) is:
\begin{equation}
\label{rfQ} \nu_l(Q^n)=2^{n-l}\binom{n}{l},\quad
l=0,\,1,\,\dots\,,\,n.
\end{equation}
The \(l\) dimensional volume of the \(l\)-face \(v_l(Q^n)\) is:
\begin{equation}
\label{vlQ}
 v_l(Q^n)=(2\rho)^l,\quad l=0,\,1,\,\dots\,,\,n.
\end{equation}
The external angle at \(l\)-face \(\gamma_l(Q^n)\) is:
\begin{equation}
\label{eaQ} \gamma_l(Q^n)=2^{-(n-l)},\quad l=0,\,1,\,\dots\,,\,n.
\end{equation}
According to \eqref{IVRP},
\(V_r(Q^n)=2^{(n-r)}\binom{n}{r}2^{-(n-r)}(2\rho)^r\), or
\begin{equation}
\label{ivCQ} V_r(Q^n)=\binom{n}{r}(2\rho)^r\,.
\end{equation}
Thus, for the cube \(Q^n\), \eqref{KisQ},
\begin{equation*}
M_{Q^n}(t)=\sum\limits_{0\leq{}l\leq{}n}\kappa_l\binom{n}{n-l}(2\rho)^{n-l}t^l\,,
\end{equation*}
or
\begin{equation}
\label{SMPQ}
M_{Q^n}(t)=(2\rho)^n\sum\limits_{0\leq{}l\leq{}n}\kappa_l\,j_{\,n,\,l}%
\frac{1}{l!}(2\rho)^{-l}(nt)^l\,.
\end{equation}
The shape factor \eqref{SF} for the cube \(Q^n\), \eqref{KisB}, is
\begin{equation}
\label{SFQ} \sigma_{Q_n}=n/\rho\,.
\end{equation}
Thus, the renormalized Steiner-Minkowski polynomial for the cube
 \(Q^n\), \eqref{KisQ}, is
\begin{equation}
\label{RSFQ}
\mathscr{M}_{Q^n}(\tau)=\sum\limits_{0\leq{}l\leq{}n}\,j_{\,n,\,l}\,\left(\frac{\sqrt{\pi}}{2}\right)^l%
\frac{1}{\Gamma\left(\frac{l}{2}+1\right)\,\Gamma(l+1)}\,\,\tau^l\,,
\end{equation}
where \(j_{n,l}\) is the Jensen multiplier defined in \eqref{JF}.
Taking into account \eqref{PJM}, we pass to the limit as
\(n\to\infty\) in the expression \eqref{RSFQ}. The limiting entire
function \(\mathscr{M}_{Q^ \infty}(\tau)\) is:
\begin{equation}
\label{RSFQL}
\mathscr{M}_{Q^\infty}(\tau)=\sum\limits_{0\leq{}l,{}\infty}\left(\frac{\sqrt{\pi}}{2}\right)^l%
\frac{1}{\Gamma(\frac{l}{2}+1)\Gamma(l+1)}\,\tau^l\,.
\end{equation}

\noindent%
\textbf{3. The family of Regular Cross-Polytopes \(C^n\).}
\begin{equation}
\label{KisC}
C^n=\{x=(\xi_1,\,\xi_2,\,\dots\,,\xi_n)\in\mathbb{R}^{n}:
\sum\limits_{1\leq{l}\leq{}n}|\xi_l|\leq{}\rho\}\,,\ \ \rho>0.
\end{equation}
The total number of \(l\)-faces \(\nu_l(C^n)\) is:
\begin{equation}
\label{rfC}%
 \nu_l(C^n)=2^{l+1}\binom{n}{l+1},\quad l=0,\,1,\,\dots\,,\,n-1\,.
\end{equation}
The \(l\)-dimensional volume of the \(l\)-face \(v_l(C^n)\) is:
\begin{equation}
\label{vlC} v_l(C^n)=\rho^l\frac{\sqrt{l+1}}{l!}
\end{equation}
The external angle at \(l\)-face \(\gamma_l(C^n)\) is calculated
in \cite{BeHe} (See Lemma 2.1 there):
\begin{equation}
\label{eaC}
 \gamma_l(C^n)=\frac{2^{n-l-1}}{\pi^{(n-l)/2}}\int\limits_{0}^{\infty}
 e^{-x^2}\Bigg(\int\limits_{0}^{x/\sqrt{l+1}}e^{-y^2}dy\Bigg)^{n-l-1}dx\,.
\end{equation}
So, the \(r\)-th intrinsic volume \(V_r(C^n)\) is:
\begin{multline}
\label{ivCC}
V_r(C^n)=2^{r+1}\binom{n}{r+1}\frac{\sqrt{r+1}}{r!}\frac{1}{\sqrt{\pi}}
\int\limits_{0}^{\infty}
 e^{-x^2}\Bigg(\frac{2}{\sqrt{\pi}}
  \int\limits_{0}^{x/\sqrt{r+1}}e^{-y^2}dy\Bigg)^{n-r-1}\hspace{-1.7ex}dx
 \cdot \rho^{r},\\
\hfill\hfill r=0,\,1,\,\dots\,,\,n-1;\hfill\\
 V_n(C_n)=\frac{2^n}{n!}\rho^n\,.
\end{multline}
\begin{multline}
M_{C^n}(t)=
\frac{2^n}{n!}\rho^n+\sum\limits_{1\leq{}l\leq{}n}\kappa_l
V_{n-l}(C^n)=\\ \frac{2^n}{n!}\rho^n+
\sum\limits_{1\leq{}l\leq{}n}\kappa_l\,2^{n-l}\binom{n}{n-l+1}%
\frac{\sqrt{n-l+1}}{(n-l)!}\cdot{}I_{n,\,l}
 \cdot \rho^{n-l}\,t^l\,,
\end{multline}
where
\begin{equation}
\label{InlC}%
 I_{n,\,l}=\frac{2}{\sqrt{\pi}}\int\limits_{0}^{\infty}
 e^{-x^2}\Bigg(\frac{2}{\sqrt{\pi}}
  \int\limits_{0}^{x/\sqrt{n-l+1}}e^{-y^2}dy\Bigg)^{l-1}\hspace{-1.7ex}dx,\quad
  1\leq{}l\leq{}n\,.
\end{equation}
 The shape factor \eqref{SF} for
the regular cross-polytope \(C^n\), \eqref{KisC}, is
\begin{equation}
\label{SFC} \sigma_{C_n}= \frac{n^{3/2}}{\rho}.
\end{equation}
The renormalized Steiner-Mincowski polynomial
\(\mathscr{M}_{C^n}\) for the family of the regular
cross-polytopes \(\{C^n\}_n\), \eqref{KisC}, is:
\begin{equation}
\label{RSFC}%
 \mathscr{M}_{C^n}(\tau)=1+\sum\limits_{1\leq{}l\leq{}n}
 \Big(\frac{\sqrt{\pi}}{2}\Big)^l
 (j_{n,\,l})^2\frac{\sqrt{n}}{\sqrt{n-l+1}}
 \frac{2n^{(l-1)/2}}{\Gamma(\frac{l}{2})}\,I_{n,\,l}\,\frac{\tau^{l}}{l!}\,.
\end{equation}
To pass to the limit as \(n\to\infty\) in \eqref{RSFQ}, we need
some information about the values \(I_{n,\,l}\), \eqref{InlC}.
Since
\(\int\limits_{0}^{\infty}e^{-\lambda^2}d\lambda=\sqrt{\pi}/2\),
\begin{equation*}%
|I_{n,\,l}|<1,\quad{}l=1,\,2,\,\dots\,,\,n\,.
\end{equation*}%
Moreover, for every fixed \(l\),
\begin{equation*}%
I_{n,\,l}=\Bigg(\frac{2}{\sqrt{\pi}}\Bigg)^l\int\limits_{0}^{\infty}
 e^{-x^2}\Bigg(\frac{x}{\sqrt{n-l+1}}\Bigg)^{l-1}dx\,(1+o(1)),\ \
 \text{as}\ \ n\to\infty\,,
\end{equation*}%
or
\begin{equation*}%
I_{n,\,l}=\Bigg(\frac{2}{\sqrt{\pi}}\Bigg)^l
\frac{1}{n^{\frac{l-1}{2}}}\int\limits_{0}^{\infty}
 x^{l-1}\,e^{-x^2}\,dx\,(1+o(1)),\ \
 \text{as}\ \ n\to\infty\,,
\end{equation*}%
and finally
\begin{equation}%
\label{AsEAC}%
 I_{n,\,l}=\frac{1}{2}\Bigg(\frac{2}{\sqrt{\pi}}\Bigg)^l\frac{1}{n^{\frac{l-1}{2}}}\,%
 \Gamma\Big(\frac{l}{2}\Big)(1+o(1)),\ \
 \text{as}\ \ n\to\infty\,,
\end{equation}%
Taking into account \eqref{PJM} and \eqref{AsEAC}, we pass to the
limit as \(n\to\infty\) in \eqref{RSFC}\,. The limiting entire
function is:
\begin{equation}
\label{RSFCL}
 \mathscr{M}_{C^\infty}(\tau)=\sum\limits_{0\leq{}l<\infty}\frac{1}{\Gamma(l+1)}\,\tau^l\,.
\end{equation}
or
\begin{equation}
\label{RMPCI}%
\mathscr{M}_{C^\infty}(\tau)=\exp\{\tau\}\,.
\end{equation}

\noindent%
\textbf{4. The family of Regular Symplexes \(S^n\).}
\begin{multline}
\label{KisS}
S^n=\{x=(\xi_1,\,\xi_2,\,\dots\,,\xi_{n+1})\in\mathbb{R}^{n+1}: \
\sum\limits_{1\leq{l}\leq{}n+1}\xi_l=\rho\,,
\\
\text{and} \ \xi_l\geq 0\,\  \text{for every} \ l \},\ \ \
\end{multline}
The symplex \(S^n\) is the \(n\)-dimensional polytop, which is
considered as a subset of the \(n\)-dimensional space
\(\mathbb{R}^n\), :
\begin{equation}
\label{NStRn}
R^n=\{x=(\xi_1,\,\xi_2,\,\dots\,,\xi_{n+1})\in\mathbb{R}^{n+1}:
\sum\limits_{1\leq{l}\leq{}n+1}\xi_l=\rho\} .
\end{equation}

 The total number of \(l\)-faces \(\nu_l(S^n)\) is:
\begin{equation}
\label{rfS}%
 \nu_l(S^n)=\binom{n+1}{n-l},\quad l=0,\,1,\,\dots\,,\,n\,.
\end{equation}
The \(l\)-dimensional volume of the \(l\)-face \(v_l(S^n)\) is:
\begin{equation}
\label{vlS} v_l(S^n)=\rho^l\frac{\sqrt{l+1}}{l!}
\end{equation}
The external angle at \(l\)-face \(\gamma_l(C^n)\) is calculated
in \cite{BoHe} (See Chapter 6, Section 6.5, Theorem 3, the formula
8 on the page 283 in  \cite{BoHe}):
\begin{equation}
\label{eaS}
 \gamma_l(S^n)=\frac{1}{\sqrt{\pi}}\int\limits_{-\infty}^{\infty}
 e^{-x^2}\Bigg(\frac{1}{\sqrt{\pi}}\int\limits_{-\infty}^{x/\sqrt{l+1}}e^{-y^2}dy\Bigg)^{n-l}dx\,.
\end{equation}
Remark that
\begin{equation}
\label{SVEAS} \gamma_0(S^n)=\frac{1}{n+1},\ \
\gamma_{n-1}(S^n)=\frac{1}{2},\ \ \gamma_{n}(S^n)=1\,.
\end{equation}
 So, the \(r\)-th intrinsic volume \(V_r(S^n)\) is:
\begin{multline}
\label{ivSC}
V_r(S^n)=\binom{n+1}{r+1}\frac{\sqrt{r+1}}{r!}\frac{1}{\sqrt{\pi}}
\int\limits_{-\infty}^{\infty}
 e^{-x^2}\Bigg(\frac{1}{\sqrt{\pi}}
  \int\limits_{-\infty}^{x/\sqrt{r+1}}e^{-y^2}dy\Bigg)^{n-r}\hspace{-1.7ex}dx
 \cdot \rho^{r},\\
\hfill r=0,\,1,\,\dots\,,\,n\,.
\end{multline}
\begin{multline}
M_{S^n}(t)=\sum\limits_{0\leq{}l\leq{}n}\kappa_l V_{n-l}(S^n)=\\
=\sum\limits_{0\leq{}l\leq{}n}\kappa_l\,\binom{n+1}{n-l+1}%
\frac{\sqrt{n-l+1}}{(n-l)!}\cdot{}I_{n,\,l}
 \cdot \rho^{n-l}\,t^l\,,
\end{multline}
where
\begin{equation}
\label{InlS}%
 I_{n,\,l}=\frac{1}{\sqrt{\pi}}\int\limits_{-\infty}^{\infty}
 e^{-x^2}\Bigg(\frac{1}{\sqrt{\pi}}
  \int\limits_{-\infty}^{x/\sqrt{n-l+1}}e^{-y^2}dy\Bigg)^{l}\hspace{-0.2ex}dx,\quad
  0\leq{}l\leq{}n\,.
\end{equation}
Remark that
\begin{equation}
\label{SVEASol}%
 I_{n,\,0}=1\,,\ I_{n,\,1}=\frac{1}{2}\,, \
I_{n,\,n}=\frac{1}{n+1}\,.
\end{equation}
(These are the relations \eqref{SVEAS} in other notation).
 The shape factor \eqref{SF} for
the regular symplex \(S^n\), \eqref{KisS}, is
\begin{equation}
\label{SFS}%
 \sigma_{S^n}= \frac{n^{3/2}(n+1)^{1/2}}{\rho}.
\end{equation}
The renormalized Steiner-Minkowski polynomial
\(\mathscr{M}_{S^n}\) for the family of the regular symplexes
\(\{S^n\}_n\), \eqref{KisS}, is:
\begin{equation}
\label{RSFS}%
 \mathscr{M}_{S^n}(\tau)=\sum\limits_{0\leq{}l\leq{}n}\kappa_l
 \frac{\sqrt{n-l+1}}{\sqrt{n+1}}\,j_{n+1,\,l}\,
 j_{n,\,l}\,\frac{(n+1)^{l/2}}{n^{l/2}}
 \,I_{n,\,l}\,\frac{\tau^{l}}{l!}\,,
\end{equation}
or
\begin{equation}
\label{RSFSoe}%
 \mathscr{M}_{S^n}(\tau)=\sum\limits_{0\leq{}l\leq{}n}\pi^{l/2}\,I_{n,\,l}\,%
 \frac{\sqrt{n-l+1}}{\sqrt{n+1}}\frac{(n+1)^{l/2}}{n^{l/2}}\,j_{n+1,\,l}\,
 j_{n,\,l}\,\frac{1}{\Gamma(\frac{l}{2}+1)\Gamma(l+1)}
 \tau^{l}\,,
\end{equation}
 To pass to the limit as \(n\to\infty\) in \eqref{RSFSoe}, we need
some information about the values \(I_{n,\,l}\), \eqref{InlS}.
Since
\(\int\limits_{-\infty}^{\infty}e^{-\lambda^2}d\lambda=\sqrt{\pi}\),
\begin{equation*}%
0<{}I_{n,\,l}<1,\quad{}l=0,\,1,\,\dots\,,\,n\,.
\end{equation*}%
Moreover, for every fixed \(l\),
\begin{equation}%
\label{AsEAS} %
I_{n,\,l}=2^{-l}(1+o(1)),\ \
 \text{as}\ \ n\to\infty\,.
\end{equation}%
Taking into account \eqref{PJM} and \eqref{AsEAS}, we pass to the
limit as \(n\to\infty\) in \eqref{RSFS}\,. The limiting entire
function is:
\begin{equation}
\label{RSFSL}
 \mathscr{M}_{S^\infty}(\tau)=
 \sum\limits_{0\leq{}l<\infty} \Big(\frac{\sqrt{\pi}}{2}\Big)^l
 \frac{\tau^l}{\Gamma(\frac{l}{2}+1)\Gamma(l+1)}\,.
\end{equation}

\section{Discussion.}
\setcounter{equation}{0} %
We have consider four families \(\{K^n\}_n\) of convex sets:
\begin{itemize}
\item[\(\circ\)]
balls \(\{B^n\}_n\), \eqref{KisB};
\item[\(\circ\)]
cubes \(\{Q^n\}_n\), \eqref{KisQ};
\item[\(\circ\)]
regular cross-polytopes \(\{C^n\}_n\), \eqref{KisC};
\item[\(\circ\)]
regular symplexes \(\{S^n\}_n\), \eqref{KisS}\,.
\end{itemize}
For every of these families and for every convex set \(K_n\) of
these family, we calculated the Steiner-Minkowski polynomial
\(M_{K_n}(t)\) explicitly. The expression of \(M_{K_n}(t)\)
contains the parameter \(\rho\), which can be considered as an
inherent linear size of \(K_n\). For example, for the ball
\(B_n\), \eqref{KisB}, \(\rho\) is the radius of the ball; for the
cube \(Q_n\), \eqref{KisQ}, \(2\rho\) is the length of its edge,
for the regular cross-polytope \(C^n\), \(\sqrt{2}\rho\) is the
length of the edge. The choice of \(\rho\) is rather uncertain. If
\(\rho\) is replaced with \(a\rho\), the value \(\rho^n\), which
appears in the expression for \(\Vol_n(K^n)\), will be replaced
with \(a^n\rho^n\), and the shape factor \(\sigma_{K^n}\) will be
replaced with \(a^{-1}\sigma_{K^n}\). We are interested in how the
values \(\Vol_n(K^n)\) and \(\sigma_{K^n}\) depend on \(n\). In
view of aforesaid, the factors of the form \(a^n\) in the
expression for \(\Vol_n(K^n)\) and \(a\) in the expression for
\(\sigma_{K^n}\) may be considered as non-essential. This
uncertainty in the choice of \(\rho\) disappear after the
renormalization of Minkowski polynomials. \textit{The renormalized
Minkowski polynomial \(\mathscr{M}_{K^n}(\tau)\) describe somehow
the shape of \(K^n\), and do not depend on the size of \(K^n\).}
The limiting entire function \(\mathscr{M}_{K^\infty}(\tau)\)
describes somehow the shape of the convex set \(K^n\) of the
family \(\{K^n\}_n\) in the very high dimension.

It looks very unexpectedly that \textit{the limiting entire
functions for the families of balls and cross-polytopes coincide}:
\begin{equation}%
\label{LEF1}%
\mathscr{M}_{B^\infty}(\tau)=\mathscr{M}_{C^\infty}(\tau)=
\sum\limits_{0\leq{}l<\infty}\frac{1}{\Gamma(l+1)}\,\tau^l,
 \end{equation}%
and that \textit{the limiting entire functions for the families of
cubes and symplexes coincide}:
\begin{equation}%
\label{LEF2}%
\mathscr{M}_{Q^\infty}(\tau)=\mathscr{M}_{S^\infty}(\tau)=
\sum\limits_{0\leq{}l<\infty}\Big(\frac{\sqrt{\pi}}{2}
\Big)^l\frac{1}{\Gamma(\frac{l}{2}+1)\Gamma(l+1)}\,\tau^l\,.
\end{equation}
Recall that the appropriate shape factors are:
\begin{equation}
\label{MESF}%
 \sigma_{B^n}=n,\ \ \sigma_{Q^n}=n,\ \
\sigma_{C^n}=n^{3/2}, \ \ \sigma_{S^n}=n^{3/2}(n+1)^{1/2}\,.
\end{equation}
(We omitted multiplicative constants, and retain the dependence on
\(n\) only in the expressions \eqref{MESF}.)

Now we focus our attention on the location of zeros of Minkowski
polynomials for the considered families of regular polytopes. Both
entire functions which appear in the right hand sides of
\eqref{LEF1} and \eqref{LEF2} belongs to the Laguerre-Polya class
\(\mathscr{L}\text{-}\mathscr{P}\text{-}\textup{I}\) of entire
functions. (See \cite{RaSc} for the definition. See also numerous
papers by Th.\,Craven and G.\,Csordas.b) For the function
\[E_1(\tau)=\sum\limits_{0\leq{}l<\infty}\frac{1}{\Gamma(l+1)}\,\tau^l\,,\]
 this is evident,
for the function
\[E_2(\tau)=\sum\limits_{0\leq{}l<\infty}\Psi(l)\frac{1}{\Gamma(l+1)}\,\tau^l\,,\]
where
\[\Psi(\tau)=\Big(\frac{\sqrt{\pi}}{2}
\Big)^\tau\,\frac{1}{\Gamma(\frac{\tau}{2})+1},\] this is a
consequence of the Laguerre Theorem, since the function
\(\Psi(\tau)\) belongs to the Laguerre-Polya class
\(\mathscr{L}\text{-}\mathscr{P}\text{-}\textup{I}\,.\) Since for
the families of balls and cubes, the renormalized Minkowski
polynomials \(\mathscr{M}_{K^n}(\tau)\) can be obtained from the
entire function
\begin{equation}%
\label{JR1}%
\mathscr{M}_{K^\infty}(\tau)=\sum\limits_{0\leq{}l<infty}\frac{\mu_l}{l!}\,\tau^l
\end{equation}%
of the class \(\mathscr{L}\text{-}\mathscr{P}\text{-}\textup{I}\)
according to the rule
\begin{equation}
\label{JR2}%
\mathscr{M}_{K^n}(\tau)=\sum\limits_{0\leq{}l<\infty}j_{n,\,l}\frac{\mu_l}{l!}\,\tau^l\,,
\end{equation}
where \(j_{n,\,l}\) are Jensen multipliers, \eqref{JF}, the zeros
of the polynomials \(\mathscr{M}_{B^n}(\tau)\) and
\(\mathscr{M}_{Q^n}(\tau)\) are simple and negative. (For the
polynomials \(\mathscr{M}_{B^n}(\tau)\), \eqref{RMPB}, this is
evident.) However, the polynomials \(\mathscr{M}_{C^n}(\tau)\) and
\(\mathscr{M}_{S^n}(\tau)\) are obtained from the appropriate
entire functions \(\mathscr{M}_{C^\infty}(\tau)\) and
 \(\mathscr{M}_{S^\infty}(\tau)\) (of the class of the class
 \(\mathscr{L}\text{-}\mathscr{P}\text{-}\textup{I}\)) by the rule more complicated  than
\eqref{JR1}\,-\,\eqref{JR2}. So, we can not conclude in the above
described way that the zeros of the polynomials
\(\mathscr{M}_{C^n}(\tau)\) and  \(\mathscr{M}_{S^n}(\tau)\) are
negative.\\

\noindent \textsf{OPEN PROBLEMS.} \\
\textbf{1.} \textit{Are all zeros of polynomials
 \(\mathscr{M}_{C^n}(\tau)\) and  \(\mathscr{M}_{S^n}(\tau)\)
negative?}\\
\textbf{2.} \textit{Are all zeros of these polynomials located in
the left half plane?}

\renewcommand{\refname}{\hfill\textsf{BIBLIOGRAPHY}\hfill\hfill}


\begin{thebibliography}{MoMo}
\bibitem[BeHe]{BeHe}\textsc{Betke, U., Henk, M.}
\textit{Intrinsic volumes and lattice points of crosspolytopes.}
Monatshefte f\"ur Math., \textbf{115}, (1993), p. 27\,--\,33.
\bibitem[B\"oHe]{BoHe}\textsc{B\"ohm, J., Hertel,\,E.},
\textit{Polyedergeometrie.} (In German).\hspace{6.0ex}\\ %
 Deutscher Verlag der
Wissenschaften, Berlin, 1980. 301 pp.
\bibitem[Gr]{Gr}\textsc{Gr\"unbaum, B.},
\textit{Convex Polytopes}. Springer-Verlag, New York, 2003.
xvi+468 pp.
\bibitem[McM]{McM}\textsc{McMullen, P.}
\textit{Non-linear angle-sum relations for polyhedral cones and polytopes.}
Math.Proc. Camb. Phil. Soc., \textbf{78} (1975), p. 247\,--\,251.
\bibitem[Ru]{Ru}\textsc{Ruben, H.} \
\textit{On the geometrical moments of skew-regular simplices in hyperspherical
space, with some applications in geometry and mathematical statistics.} Acta Math.,
\textbf{103} (1960), p.\,1\,--\,23.
\bibitem[RaSc]{RaSc}\textsc{Rahman,\,Q.L.} and
\textsc{G,\,Schmeisser.} \textit{Analytic Theory of Polynomials.}
Clarendon Press, Oxford, 2002. xiv+742.
\bibitem[Schn]{Schn}\textsc{Schneider, R.}\ \  \textit{Convex Bodies:
The Brunn-Minkowski Theory.} Cambridge Univ.
Press, Cambridge 1993.
\end{thebibliography}
\end{document}